\documentclass[a4paper, 11pt, twoside, final]{scrartcl}

\usepackage{PREAMBLE}
\usepackage{COMMAND}
\usepackage{HYPHENATION}

\begin{document}
\thispagestyle{empty}
\begin{center}
  \rule{\linewidth}{1pt}\\[0.4cm]
  {\sffamily \bfseries \large Sparse phase retrieval of one-dimensional signals 
    by Prony's method}\\[10pt] 
  {\sffamily Robert Beinert\textsuperscript{1} and
    Gerlind Plonka\textsuperscript{2}}\\[5pt]
  {\sffamily\scriptsize \textsuperscript{1}Institut für Mathematik und
  Wissenschaftliches Rechnen}\\[-3pt] 
  {\sffamily\scriptsize Karl-Franzens-Universität Graz}\\
  {\sffamily\scriptsize \textsuperscript{2}Institut für Numerische und
    Angewandte Mathematik}\\[-3pt] 
  {\sffamily\scriptsize Georg-August-Universität Göttingen}\\
  \rule{\linewidth}{1pt}
\end{center}

\vspace*{10pt}

{\small
  \noindent
  {\sffamily\bfseries Abstract:} In this paper, we show that sparse
  signals $f$ representable as a linear combination of a finite number
  $N$ of spikes at arbitrary real locations or as a finite linear
  combination of B-splines of order $m$ with arbitrary real knots can
  be almost surely recovered from ${\cal O}(N^{2})$ intensity
  measurements
  $\abs{\Fourier \mleft[ f \mright] \mleft( \omega \mright) }^{2}$ up
  to trivial ambiguities. The constructive proof consists of two
  steps, where in the first step the Prony method is applied to
  recover all parameters of the autocorrelation function and in the
  second step the parameters of $f$ are derived.  Moreover, we present
  an algorithm to evaluate $f$ from its Fourier intensities and
  illustrate it at different numerical examples.

  \smallskip


  \smallskip
  
}

\section{Introduction}
\label{sec:introduction}
\setcounter{equation}{0}

Phase retrieval problems occur in many scientific fields, particularly
in optics and communications. They have a long history with rich
literature regarding uniqueness of solutions and existence of reliable
algorithms for signal reconstruction, see e.g. \cite{SECCMS15} and
references therein.  Usually, the challenge in solving one-dimensional
phase retrieval problems is to overcome the strong ambiguousness by
determining appropriate further information on the solution signal.
Previous literature on characterization of ambiguities of the phase
retrieval problem with given Fourier intensities is often concerned
with the discrete problem, where a signal ${\bf x}$ in ${\R}^{N}$ or
${\C}^{N}$ has to be recovered. For an overview on the complete
characterization of nontrivial ambiguities is this discrete case as
well as on appropriate additional signal information we refer to our
survey \cite{BP15} and further recent results in
\cite{BP16,Bei16-1,Bei16}.

\noindent
{\bf Contribution of this paper.}
In this paper, we consider the continuous one-dimensional sparse phase
retrieval problem to reconstruct a complex-valued signal from the
modulus of its Fourier transform.  Applications of this problem occur
in electron microscopy, wave front sensing, laser optics
\cite{SST04,SSD+06} as well as in X-ray crystallography and speckle
imaging \cite{RCLV13}.  For the posed problem, we will show that for
sparse signals the given Fourier intensities are already sufficient
for an almost sure unique recovery, and we will give a construction algorithm to
recover $f$.

We assume that the sparse signal is either of the form
\begin{equation} \label{spikes}
 f(t) = \sum_{j=1}^{N} c_{j}^{(0)} \, \delta( t - T_{j}) 
\end{equation} 
or, for $m >0$,
\begin{equation} \label{splines}
 f(t) = \sum_{j=1}^{N} c_{j}^{(m)} \, B_{j,m}(t) 
 \end{equation}
 with $c_{j}^{(m)} \in \C$, $T_{j} \in {\R}$ for $j=1, \ldots , N$,
 where $\delta$ denotes the Delta distribution, and $B_{j,m}$ is the
 B-spline of order $m$ being determined by the (real) knots
 $T_{j} < T_{j+1} < \ldots < T_{j+m}$.  We want to recover these
 signals from the Fourier intensities $| \widehat{f}(\omega)|^{2}$ and
 will show that only ${\cal O}(N^{2})$ samples are needed to recover
 $f$, i.e. all coefficients $c_{j}^{(m)}$, $j=1, \ldots, N$ and knots
 $T_{j}$, $j=1, \ldots , N+m$, almost surely up to trivial
 ambiguities.  The proposed procedure is constructive and consists in
 two steps. In a first step, we employ Prony's method in oder to
 recover all parameters of the (squared) Fourier intensity function
 $\abs{\Fourier \mleft[ f \mright] \mleft( \omega \mright) }^{2}$. In
 a second step, we recover the parameters $T_{j}$ and the complex
 coefficients $c_{j}$ that determine the desired signal.


\noindent
{\bf Related work on sparse phase retrieval.}
While the general phase retrieval problem has been extensively studied
for a long tome, the special case of sparse phase retrieval grew to a
strongly emerging field of research only recently, particularly often
connected with ideas from compressed sensing.  Most of the papers
consider a discrete setting, where the $N$-dimensional real or complex
$k$-sparse vector ${\bf x}$ has to be reconstructed from measurements
of the more general form
$| \langle {\bf a}_{j}, {\bf x} \rangle |^{2}$ with vectors
${\bf a}_{j}$ forming the rows of a measurement matrix
${\bf A} \in {\C}^{M \times N}$.  The needed number $M$ of
measurements depends on the sparsity $k$.

If ${\bf A}$ presents rows of a Fourier matrix, this setting is close
to the sparse phase retrieval problem considered in optics, see e.g.\
\cite{JOH13}. Here the problem is first rewritten as (non-convex) rank
minimization problem, then a tight convex relaxation is applied and
the optimization problem is solved by a re-weighted $l_{1}$-minimization
method.  The related approach in \cite{ESMBC15} employs the magnitudes
of the short-time Fourier transform and applies the occurring
redundancy for unique recovery of the desired signal. A corresponding
reconstruction algorithm is then based on an adaptation of the GESPAR
algorithm in \cite{SBE14}.

In \cite{LV13}, the measurement matrix ${\bf A}$ is taken with random
rows and the PhaseLift approach \cite{CSV11} leads to a convex
optimization problem that recovers the sparse solution with high
probability.  Employing a thresholded gradient descent algorithm to a
non-convex empirical risk minimization problem that is derived from
the phase retrieval problem, Cai et al. \cite{CLM16} have established
the minimax optimal rates of convergence for noisy sparse phase
retrieval under sub-exponential noise.

Other papers rely on the compressed sensing approach to construct
special frame vectors ${\bf a}_{j}$ to ensure uniqueness of the phase
retrieval problem with high probability, where the number of needed
vectors is ${\cal O}(k)$, see e.g. \cite{WX14,OE14,IVW17}.

We would like to emphasize that all approaches employing general or
random measurement matrices in phase retrieval are quite different in
nature from our phase retrieval problem based on Fourier intensity
measurements.  In this paper, we want to stick on considering Fourier
intensity measurements because of their particular relevance in
practice.

Early attempts to exploit sparsity of a discrete signal for unique
recovery using Fourier intensities go back to unpublished manuscripts
by Yagle \cite{Yag1,Yag2}, where a variation of Prony's method is
applied in a non-iterative algorithm to sparse signal and image
reconstruction. Unfortunately, the algorithm proposed there not always
determines the signal support correctly.

The continuous one-dimensional phase retrieval problem has been rarely
discussed in the literature, see
\cite{Wal63,Hof64,RCLV13,Bei16,BP15a}.  In the preprint \cite{RCLV13},
the authors also considered the recovery of sparse continuous signals
of the form (\ref{spikes}).  However, in that paper the sparse phase
retrieval problem is in turn transferred into a turnpike problem that
is computationally expensive to solve. Moreover there exist cases,
where a unique solution cannot be found, see \cite{Blo75}.  Our method
circumvents this problem by proposing an iterative procedure to fix
the signal support (resp. the knots of the signal represented as a
B-spline function) where the corresponding signal coefficients are
evaluated simultaneously.
 

\noindent 
{\bf Organization of this paper.}
In Section 2, we shortly recall the mathematical formulation of the
considered sparse phase retrieval problem and the notion of trivial
ambiguities of the phase retrieval problem that always occur.
 
Section 3 is devoted to the special case of phase retrieval for
signals of the form (\ref{spikes}). Using Prony's method, we give a
constructive proof for the unique recovery of the $N$-sparse signal
$f$ up to trivial ambiguities using $\nicefrac{3}{2} \, N(N-1)+1$ Fourier
intensity measurements. Here we have to assume that the knot
differences $T_{j} - T_{k}$ are pairwise different.

In Section 4, the ansatz is generalized to the unique recovery of
spline functions of the form (\ref{splines}) where we need to employ
$\nicefrac{3}{2}(N+m)(N+m-1)+1$ Fourier intensity measurements. In Section
5, we present an explicit algorithm for the considered sparse phase
retrieval problem and illustrate it at different examples.

\section{Trivial ambiguities of the phase retrieval problem}
\label{sec:phase-retrieval}
\setcounter{equation}{0}

We wish to recover an unknown complex-valued signal
$f \colon \R \rightarrow \C$ of the form (\ref{spikes}) or
(\ref{splines}) with compact support from its Fourier intensity
$\absn{\Fourier [f]}$ given by
\begin{equation*}
  \abs{\Fourier \mleft[ f \mright] \mleft( \omega \mright) }
  \coloneqq \abs{\widehat f \mleft( \omega \mright)}
  \coloneqq \absB{\int_{-\infty}^\infty f\mleft(t\mright) \, \e^{- \I
      \omega t} \diff{t}}
  \qquad
  (\omega \in \R).
\end{equation*}
Unfortunately, the recovery of the signal $f$ is complicated because of the
well-known ambiguousness of the phase retrieval problem.  Transferring
\cite[Proposition~2.1]{BP15} to our setting, we can recover $f$ only
up to the following ambiguities.

\begin{Proposition}\label{prop:trivial-amb}
  Let $f$ be of a signal of the form \eqref{spikes} or a non-uniform spline function of the form \eqref{splines}. Then
  \begin{enumerate}[\upshape(i)]
  \item\label{item:1} the rotated signal $\e^{\I \alpha} \, f$ for
    $\alpha \in \R$,
  \item\label{item:2} the time shifted signal $f( \cdot-t_0)$ for
    $t_0 \in \R$,
  \item\label{item:3} and the conjugated and reflected signal
    $\overline{ f(-\cdot)}$
  \end{enumerate}
  have the same Fourier intensity $\absn{\Fourier [f]}$.%
\end{Proposition}

\begin{Proof}
  Applying the properties of the Fourier transform, we have
  \begin{enumerate}[(i)]
  \item
    $\Fourier[\e^{\I \alpha} \, f] = \e^{\I \alpha} \, \Fourier [f];$
  \item
    $\Fourier[f(\cdot - t_0)] = \e^{-\I \omega t_0} \, \Fourier [f];$
  \item $\Fourier[\overline{f[ -\cdot ]}] = \overline{\Fourier [f]}.$
  \end{enumerate}
  Considering the absolute value of each equation yields the
  assertion.  \qed
\end{Proof}

Although the ambiguities in Proposition~\ref{prop:trivial-amb} always
occur, they are of minor interest because of their close relation to
the original signal.  For this reason, we call ambiguities caused by
rotation, time shift, conjugation and reflection, or by combinations
of these \emph{trivial}.  In the following, we will show that for the
considered sparse signals the remaining non-trivial ambiguities only
occur in rare cases.

\section{Phase retrieval for distributions with discrete support}
\label{sec:retr-step-funct}
\setcounter{equation}{0}

Initially, we restrict ourselves to the recovery of signals $f$ of the form (\ref{spikes})
with complex-valued coefficients $c_j^{(0)}$ and spike locations
 $T_1 < \dots < T_{N}$.  
\begin{equation*}
 \widehat f \mleft(\omega\mright)
  = \sum_{j=1}^{N} c_j^{(0)} \, \e^{-\I \omega T_j}
  \qquad (\omega \in \R),
\end{equation*}
and the  known squared Fourier intensity $\absn{\Fourier [f]}^{2}$ can be represented by
\begin{equation}
  \label{eq:Four-int-spikes}
  \abs{\widehat f \mleft(\omega\mright)}^2 
  = \sum_{j=1}^{N} \sum_{k=1}^{N} c_j^{(0)} \overline c_k^{(0)} \, \e^{-\I \omega (T_j - T_k)}.
\end{equation}
Thus, in order to recover $f$ being determined by the coefficients
$c_{j}^{(0)} \in {\C}$ and the knots $T_{j} \in {\R}$,
$j=1, \ldots , N$, we will recover all parameters of the exponential
sum in (\ref{eq:Four-int-spikes}) in a first step and then derive the
desired parameters of $f$ in a second step.

\subsection{First step: Parameter recovery by Prony's method}
\label{sec:pronys-method}

Assuming that the non-zero knot differences $T_j - T_k$ with
  $j \ne k$ are pairwise different, and denoting the distinct
  frequencies $T_j - T_k$ in increasing order by $\tau_\ell$ with
  $\ell = - \nicefrac{N(N-1)}{2}, \dots, \nicefrac{N(N-1)}{2}$, we can
  rewrite (\ref{eq:Four-int-spikes}) as
\begin{equation}\label{ex1}
 P(\omega) := \abs{\widehat f \mleft(\omega\mright)}^2 
  = \sum_{\ell=-\nicefrac{N(N-1)}{2}}^{\nicefrac{N(N-1)}{2}}  \gamma_{\ell} \, \e^{-\I \omega \tau_{\ell}}
  = \gamma_{0} + \sum_{\ell=1}^{\nicefrac{N(N-1)}{2}}  \left( \gamma_{\ell} \, \e^{-\I \omega \tau_{\ell}} + \overline{\gamma}_{\ell} \, \e^{\I \omega \tau_{\ell}} \right)
\end{equation}
with the related coefficients
$\gamma_\ell \coloneqq c_j^{(0)} \overline c_k^{(0)}$ for the non-zero
frequencies $\tau_\ell = T_j - T_k$ and
\raisebox{0pt}[0pt][0pt]{$\gamma_0 \coloneqq \sum_{j=1}^N
  \absn{c_j^{(0)}}^2$}
for the zero frequency.  Since $\tau_{-\ell} = - \tau_\ell$, the
coefficients in \eqref{ex1} fulfill the conjugated symmetry
$\gamma_{-\ell} = \overline \gamma_\ell$.

In order to recover the parameters $\tau_{\ell}$ and the unknown coefficients $\gamma_{\ell}$ 
 from the exponential sum \eqref{ex1} we employ 
Prony's method \cite{Hil87, PT14}. 
Let $h>0$ be chosen such that $h\tau_{\ell} <  \pi$ for all $\ell=1, \ldots , \nicefrac{N(N-1)}{2}$.

Using the intensity values
$P(hk) = \abs{\widehat f \mleft(h k \mright)}^2 $, $k = 0, \dots, 2N(N-1)+1$, the unknown
parameters $\gamma_{\ell}$ and $\tau_{\ell}$ in (\ref{ex1}) can be determined by exploiting the
algebraic Prony polynomial $\Lambda(z)$ defined by
\begin{equation}
  \label{eq:Prony-poly}
  \Lambda \mleft( z \mright)
  \coloneqq
  \prod_{\ell=-\nicefrac{N(N-1)}{2}}^{\nicefrac{N(N-1)}{2}} \left( z - \e^{-\I h \tau_\ell} \right)
  = \sum_{k=0}^{N(N-1)+1} \lambda_k \, z^k,
\end{equation}
where $\lambda_k$ denote the coefficients in the monomial representation of  $\Lambda(z)$.
Obviously, $\Lambda(z)$ is always a monic polynomial, which means that
$\lambda_{N(N-1)+1} = 1$.

Using the definition of the Prony polynomial $\Lambda(z)$ in (\ref{eq:Prony-poly}), we observe
that 
\begin{eqnarray*}
  \sum_{k=0}^{N(N-1)+1} \lambda_k \, P\mleft( h \mleft( k + m \mright) \mright)
  &=& \sum_{k=0}^{N(N-1)+1} \sum_{\ell=-\nicefrac{N(N-1)}{2}}^{\nicefrac{N(N-1)}{2}} \lambda_k \gamma_{\ell} \, \e^{-\I h (k+m) \tau_{\ell}} \\
  &=& \sum_{\ell=-\nicefrac{N(N-1)}{2}}^{\nicefrac{N(N-1)}{2}} \gamma_{\ell} \, \e^{-\I h m \tau_{\ell}} \,
  \Lambda \mleft( \e^{- \I h \tau_{\ell}}
  \mright)
  =0
\end{eqnarray*}
for $m=0, \dots, N(N-1)$.  Consequently, the vector of remaining coefficients
$\Vek \lambda \coloneqq (\lambda_0, \dots, \lambda_{N(N-1)})^\T$ of the
Prony polynomial $\Lambda(z)$ can be determined by solving the linear equation
system
\begin{equation}
  \label{eq:Hankel-system}
  \Mat H  \Vek \lambda = -\Vek h
\end{equation}
with $\Mat H \coloneqq (P(h(k+m)))_{m,k=0}^{N(N-1)}$ and
$\Vek h \coloneqq (P(h(N(N-1)+1+m)))_{m=0}^{N(N-1)}$.  Since the Hankel matrix
$\Mat H$ can be written as
\begin{equation*}
  \Mat H = \Mat V^\T \diag \mleft( \gamma_{-\nicefrac{N(N-1)}{2}}, \dots, \gamma_{\nicefrac{N(N-1)}{2}} \mright) \, \Mat V
  \addmathskip
\end{equation*}
with the Vandemonde matrix
$\Mat V \coloneqq (\e^{-h k
  \tau_{\ell}})_{\ell=-\nicefrac{N(N-1)}{2},k=0}^{\nicefrac{N(N-1)}{2},N(N-1)+1}$,
the linear equation system \eqref{eq:Hankel-system} possesses a unique
solution if and only if the unimodular values $\e^{-\I h \tau_{\ell}}$
differ pairwise for
$\ell=-\nicefrac{N(N-1)}{2}, \dots, \nicefrac{N(N-1)}{2}$.  This
assumption has been ensured by choosing an $h$ such that
$h \tau_{\ell} \in (- \pi, \pi)$, since the $\tau_{\ell}$ had been
supposed to be pairwise different.

Knowing the coefficients $\lambda_k$ of $\Lambda(z)$, we can determine
the unknown frequencies $\tau_{\ell}$ by evaluating the roots of the Prony polynomial
\eqref{eq:Prony-poly}.  The coefficients $\gamma_{\ell}$ can now be
computed by solving the over-determined equation system
\begin{equation} \label{van}
  \sum_{\ell=-\nicefrac{N(N-1)}{2}}^{\nicefrac{N(N-1)}{2}}  \gamma_{\ell} \, \e^{-\I h k \tau_{\ell}}
  = P \mleft( h k \mright)
  \qquad
  (k=0, \dots, 2N(N-1)+1)
\end{equation}
with a Vandermonde-type system matrix.

The procedure summarized above is the usual Prony method, adapted to the non-negative exponential sum $P(\omega)$ in (\ref{ex1}).
In the numerical experiments in Section 5, we will apply the approximate Prony method (APM) in \cite{PT10}. APM is based on the above considerations but it is numerically more stable 
and exploits the special properties $\gamma_{-\ell} = \overline{\gamma}_{\ell}$ and $\tau_{-\ell} = - \tau_{\ell}$ for $\ell=0, \ldots , \nicefrac{N(N-1)}{2}$.

Let us now investigate the question, how many intensity values are at
least necessary for the recovery of $P(\omega)$ in (\ref{ex1}).
Counting the number of unknowns of $P(\omega)$ in (\ref{ex1}), we only
need to recover the $\nicefrac{3}{2} \, N(N-1) +1$ real values
$\gamma_{0}$ and $\re \gamma_{\ell}$, $\im \gamma_{\ell}$,
$\tau_{\ell}$, for $\ell=1, \ldots \nicefrac{N(N-1)}{2}$.  We will
show now that using the special structure of the real polynomial
$P(\omega)$ in (\ref{ex1}) and of the Prony polynomial $\Lambda(z)$ in
(\ref{eq:Prony-poly}), we indeed need only
$\nicefrac{3}{2} \, N(N-1) +1$ exact equidistant real measurements
$P(kh)$, $k=0, \ldots , \nicefrac{3}{2} \, N(N-1)$ to recover all
parameters determining $P(\omega)$. This can be seen as follows.

Reconsidering $\Lambda(z)$ in (\ref{eq:Prony-poly}) with $\tau_{0}=0$
and $\tau_{\ell} = - \tau_{-\ell}$, we obtain
\begin{eqnarray*}
\Lambda \mleft( z \mright)
  &=& (z-1) \prod_{\ell=1}^{\nicefrac{N(N-1)}{2}} \left( z - \e^{\I h \tau_\ell} \right) \left( z - \e^{-\I h \tau_\ell} \right) \\
  &=& (z-1) \prod_{\ell=1}^{\nicefrac{N(N-1)}{2}} \left( z^{2} - 2z \cos(h \tau_{\ell}) +1  \right) 
  = \sum_{k=0}^{N(N-1)+1} \lambda_k \, z^k,
\end{eqnarray*}
where all occurring coefficients $\lambda_{k}$ are real.  Moreover, since
$$ z^{-\nicefrac{(N(N-1)+1)}{2}} \Lambda \mleft( z \mright)
= (z^{\nicefrac{1}{2}} - z^{-\nicefrac{1}{2}}) \prod_{\ell=1}^{\nicefrac{N(N-1)}{2}} \left( z - 2 \cos(h \tau_{\ell}) +z^{-1}  \right) 
$$
is antisymmetric, it follows that 
$$ \lambda_{N(N-1)+1-k} = - \lambda_{k} \qquad (k=0, \ldots , \nicefrac{N(N-1)}{2}),
$$
and particularly $\lambda_{N(N-1)+1} = -\lambda_{0} =1$.
In order to determine the unknown coefficients $\lambda_{k}$, $k=1, \ldots, \nicefrac{N(N-1)}{2}$ of 
$$ \Lambda \mleft( z \mright) = \sum_{k=0}^{\nicefrac{N(N-1)}{2}} \lambda_{k} \left( z^{k} - z^{N(N-1)+1-k} \right),
$$
we employ (\ref{ex1}) and observe that for $m=0, \ldots , \nicefrac{N(N-1)}{2} -1$,
\begin{eqnarray*}
& & \sum_{k=0}^{\nicefrac{N(N-1)}{2}} \lambda_{k} \left[ P(h(k+m)) - P(h(N(N-1)+1+m-k)) \right] \\
&=& \sum_{k=0}^{\nicefrac{N(N-1)}{2}}  \lambda_{k} \left[ \sum_{\ell=1}^{\nicefrac{N(N-1)}{2}} \gamma_{\ell} \left( \e^{-\I h (k+m) \tau_{\ell}} - \e^{-\I h (N(N-1)+1+m-k) \tau_{\ell}} \right) \right. \\
& & + \left. \sum_{\ell=1}^{\nicefrac{N(N-1)}{2}} \overline{\gamma}_{\ell} \left( \e^{\I h (k+m) \tau_{\ell}} - \e^{\I h (N(N-1)+1+m-k) \tau_{\ell}} \right) \right] \\
&=& \sum_{\ell=1}^{\nicefrac{N(N-1)}{2}} \gamma_{\ell} \, \e^{-\I h m \tau_{\ell}}  \sum_{k=0}^{\nicefrac{N(N-1)}{2}} \lambda_{k} \left( \e^{-\I h k \tau_{\ell}} - 
\e^{-\I h (N(N-1)+1-k)\tau_{\ell}} \right) \\
& & + \sum_{\ell=1}^{\nicefrac{N(N-1)}{2}} \overline{\gamma}_{\ell} \, \e^{\I h m \tau_{\ell}} \sum_{k=0}^{\nicefrac{N(N-1)}{2}} \lambda_{k} \left( \e^{\I h k \tau_{\ell}} - 
\e^{\I h (N(N-1)+1-k)\tau_{\ell}} \right) \\
&=&  \sum_{\ell=1}^{\nicefrac{N(N-1)}{2}} \gamma_{\ell} \, \e^{-\I h m \tau_{\ell}} \Lambda(\e^{-\I h \tau_{\ell}}) + 
\sum_{\ell=1}^{\nicefrac{N(N-1)}{2}} \overline{\gamma}_{\ell} \, \e^{\I h m \tau_{\ell}} \Lambda(\e^{\I h \tau_{\ell}}) = 0.
\end{eqnarray*}
Therefore, the vector of unknown coefficients $\Vek \lambda \coloneqq (\lambda_1, \dots, \lambda_{\nicefrac{N(N-1)}{2}})^\T$ can be already 
evaluated from the system 
\begin{align*}
    & \sum_{k=1}^{\nicefrac{N(N-1)}{2}} \lambda_{k} \left[ P(h(k+m)) -P(h(N(N-1)+1+m-k)) \right] \\
    & \qquad \qquad = \left[ P(hm) - P(h(N(N-1)+1+m)) \right]
\hspace{60pt}
(m=0, \ldots , \nicefrac{N(N-1)}{2} -1).
\end{align*}
The parameters $\tau_{\ell}$ are then extracted from the zeros of $\Lambda(z)$, and the coefficients $\gamma_{\ell}$, $\ell=0, \ldots , \nicefrac{N(N-1)}{2}$,
are computed as in (\ref{van}) but with  $k=0, \ldots ,
\nicefrac{3}{2} \, N(N-1)$.

\subsection{Second step: Unique signal recovery}

Having determined  the parameters $\tau_{\ell}$  as well as the
corresponding coefficients $\gamma_{\ell}$ of (\ref{ex1}), we want to
reconstruct the parameters $T_{j}$ and $c_{j}^{(0)}$, $j=1, \ldots , N$,  of $f$ in (\ref{spikes}) in a second step.

\begin{Theorem}
  \label{the:uni-step-fct}
  Let $f$ be a signal  of the form \eqref{spikes}, whose
  knot differences $T_j - T_k$ differ pairwise for
  $j,k \in \{1,\dots,N\}$ with $j \ne k$, and whose coefficients
  satisfy
  \raisebox{0pt}[0pt][0pt]{$\absn{c^{(0)}_1} \ne \absn{c^{(0)}_N}$}.
  Further, let $h$ be a step size such that
  $h (T_j - T_k) \in (-\pi, \pi)$ for all $j, k$.  Then $f$ can be uniquely
  recovered from its Fourier intensities $\absn{\Fourier [f](h\ell)}$
  with $\ell = 0, \dots, \nicefrac{3}{2} \, N(N-1)$ up to trivial ambiguities.
\end{Theorem}

\begin{Proof}
  Applying Prony's method to the given data
  $\absn{\Fourier [f](h \ell)}$, we can compute the frequencies
  $\tau_\ell$ and the related coefficients $\gamma_\ell$ of the
  squared Fourier intensity \eqref{ex1}.  Again, we assume that the
  frequencies $\tau_\ell$ occur in increasing order and, further,
  denote the list of positive frequencies by
  ${\cal T} \coloneqq \{ \tau_{\ell}
  \}_{\ell=1}^{\nicefrac{N(N-1)}{2}}$.

  Obviously, the maximal distance $\tau_{\nicefrac{N(N-1)}{2}}$ now corresponds to
  the length $T_{N} - T_{1}$ of the unknown $f$ in (\ref{spikes}).
  Due to the trivial shift ambiguity, we can assume without loss of
  generality that $T_1 = 0$ and $T_{N}=\tau_{\nicefrac{N(N-1)}{2}}$.
  Further, the second largest distance
  $\tau_{(\nicefrac{N(N-1)}{2})-1}$ corresponds either to
  $T_{N-1} - T_1$ or to $T_{N} - T_2$.  Due to the trivial reflection
  and conjugation ambiguity, we can assume that
  $T_{N-1} = \tau_{(\nicefrac{N(N-1)}{2})-1}$.  By definition, there
  exists a $\tau_{\ell^*} > 0$ in our sequence of parameters
  ${\cal T}$ such that
  $ \tau_{\ell^*} + \tau_{(\nicefrac{N(N-1)}{2})-1} =
  \tau_{\nicefrac{N(N-1)}{2}}$,
  and $\tau_{\ell^*}$ hence corresponds to the knot difference
  $T_N - T_{N-1}$. Thus, we obtain
  \begin{equation*}
    c^{(0)}_{N} \overline c^{(0)}_1= \gamma_{\nicefrac{N(N-1)}{2}},
    \qquad
    c^{(0)}_{N-1} \overline  c^{(0)}_1 = \gamma_{(\nicefrac{N(N-1)}{2})-1},
    \quad\text{and}\quad
    c^{(0)}_{N} \overline c^{(0)}_{N-1} = \gamma_{\ell^*}.
  \end{equation*}
   These equations lead us to
  \begin{equation*}
    c^{(0)}_{N} = \tfrac{\gamma_{\nicefrac{N(N-1)}{2}}}{\overline c^{(0)}_1},
    \qquad
    c^{(0)}_{N-1}  = \tfrac{\gamma_{(\nicefrac{N(N-1)}{2})-1}}{\overline
      c^{(0)}_1},
  \end{equation*}
  and thus to
  \begin{equation*}
    \abs{c_1^{(0)}}^2
    = \tfrac{\gamma_{\nicefrac{N(N-1)}{2}} \overline{\gamma_{(\nicefrac{N(N-1)}{2})-1}}}{\gamma_{\ell^*}}.
    \addmathskip
  \end{equation*}
  Since $f$ can only be recovered up to a global rotation, we can
  assume that $c_1^{(0)}$ is real and non-negative, which allows us to
  determine the coefficients $c_1^{(0)}$, $c_N^{(0)}$, and
  $c_{N-1}^{(0)}$ in a unique way.

  To determine the remaining coefficients and support knots of $f$, we
  notice that the third largest distance
  $\tau_{(\nicefrac{N(N-1)}{2})-2}$ corresponds either to
  $T_{N} - T_2$ or to $T_{N-2}-T_1$.  As before, we always find a
    frequency $\tau_{\ell^*}$ such that
    $\tau_{(\nicefrac{N(N-1)}{2})-2} + \tau_{\ell^*} =
    \tau_{\nicefrac{N(N-1)}{2}}$.
  
  \noindent
  Case 1: If $\tau_{(\nicefrac{N(N-1)}{2}) - 2} = T_{N} - T_{2}$, then
  we have 
  $$\tau_{\ell^*} = \tau_{\nicefrac{N(N-1)}{2}} -
  \tau_{(\nicefrac{N(N-1)}{2})-2}  = (T_{N}- T_{1}) - (T_{N} - T_{2})
  =  T_{2} - T_{1}$$
  with  the related coefficient $\gamma_{\ell^*} = c_{2}^{(0)} \overline{c}_{1}^{(0)}$. Moreover,
  we have $\gamma_{(\nicefrac{N(N-1)}{2})-2} = c_{N}^{(0)} \overline{c}_{2}^{(0)}$ such that
  \begin{equation}
    \label{eq:case1}
    c_{2}^{(0)} = \frac{\gamma_{\ell^*} }{\overline{c}_{1}^{(0)}} = \frac{\overline{\gamma}_{(\nicefrac{N(N-1)}{2})-2}}{\overline{c}_{N}^{(0)}}.
  \end{equation}
  \noindent
  Case 2: If $\tau_{(\nicefrac{N(N-1)}{2}) - 2} = T_{N-2} - T_{1}$,
  then we have
  $$\tau_{\ell^*} = \tau_{\nicefrac{N(N-1)}{2}} -
  \tau_{(\nicefrac{N(N-1)}{2})-2} = (T_{N}- T_{1}) - (T_{N-2} - T_{1}) = T_{N} - T_{N-2} $$
  with  the related coefficient $\gamma_{\ell^*} = c_{N}^{(0)} \overline{c}_{N-2}^{(0)}$ and 
  $\gamma_{(\nicefrac{N(N-1)}{2})-2} = c_{N-2}^{(0)} \overline{c}_{1}^{(0)}$. Thus,
  \begin{equation}
    \label{eq:case2}
    c_{N-2}^{(0)} = \frac{\overline{\gamma}_{\ell^*} }{\overline{c}_{N}^{(0)}} = \frac{{\gamma}_{(\nicefrac{N(N-1)}{2})-2}}{\overline{c}_{1}^{(0)}}.
  \end{equation}
  However, only one of the two equalities in (\ref{eq:case1}) and (\ref{eq:case2}) can be true, since if both were true then
  $ \gamma_{\ell^*} \overline{c}_{N}^{(0)} = \overline{c}_{1}^{(0)} \overline{\gamma}_{(\nicefrac{N(N-1)}{2}) -2}$ and 
  $ \overline{c}_{1}^{(0)} \overline{\gamma}_{\ell^*} = \overline{c}_{N}^{(0)} \gamma_{(\nicefrac{N(N-1)}{2}) -2}$ lead to
  $$ \left| \frac{c_{N}^{(0)}}{c_{1}^{(0)}} \right| = \left| \frac{\gamma_{(\nicefrac{N(N-1)}{2}) -2}}{\gamma_{\ell^*}} \right| = \left| \frac{c_{1}^{(0)}}{c_{N}^{(0)}} \right| 
  $$
  contradicting the assumption that  $|c_{N}^{{0}}| \neq |c_{1}^{{0}}|$.
Consequently, either the equation in \eqref{eq:case1} or the
  equation in \eqref{eq:case2} holds true and we can either determine $T_{2}$ with $c_{2}^{(0)}$ or $T_{N-2}$ with $c_{N-2}^{(0)}$.
  Removing all parameters $\tau_{\ell}$ from the sequence of distances  ${\cal T}$ that correspond to 
 the differences $T_j - T_k$ of the recovered knots, we can repeat this approach to find
  the remaining coefficients and knots of $f$ inductively.  
  \qed
\end{Proof}

If we identify the space of complex-valued signals of the form
(\ref{spikes}) with the real space $\R^{3N}$, the condition that two knot
differences $T_{j_1}-T_{k_1}$ and $T_{j_2}-T_{k_2}$ are equal for
fixed indices $j_1$, $j_2$, $k_1$, and $k_2$ defines a hyper plane
with Lebesgue measure zero.  An analogous observation follows for the
condition
\raisebox{0pt}[0pt][0pt]{$\absn{c_1^{(0)}} = \absn{c_N^{(0)}}$}.  The
signals excluded in Theorem~\ref{the:uni-step-fct} hence form a negligible
null set.

\begin{Corollary}
  \label{cor:uni-step-fct}
  Almost all signals $f$ in \eqref{spikes} can be uniquely
  recovered from their Fourier intensities $\abs{\Fourier [f]}$ up to
  trivial ambiguities.
\end{Corollary}

\begin{Remark}
  1. Since the proof of Theorem~\ref{the:uni-step-fct} is
  constructive, it can be used to recover an unknown signal
  (\ref{spikes}) analytically and numerically.  If the number $N$ of
  spikes is known beforehand then the assumption of
  Theorem~\ref{the:uni-step-fct} can be simply checked during the
  computation.  If the assumption regarding pairwise different
  distances $T_{j}-T_{k}$ is not satisfied, then the application of
  Prony's method in the first step yields less than $N(N-1)+1$
  pairwise distinct parameters $\tau_{\ell}$.  The second assumption
  $|c_{N}^{{0}}| \neq |c_{1}^{{0}}|$ can be verified in the second
  step, where $c_{1}^{(0)}$, $c_{N-1}^{(0)}$, and $c_{N}^{(0)}$ are
  evaluated.
  
  2. A similar phase retrieval problem had been transferred to a turnpike problem in \cite{RCLV13}.
  The turnpike problem deals  with the recovery of the knots $T_{j}$ from an unlabeled set of distances.
  Although this problem is solvable under certain conditions, a  backtracing algorithm can have exponential complexity in the worst case, see
  \cite{LSS03}. 
\end{Remark}

\section{Retrieval of spline functions with arbitrary knots}
\setcounter{equation}{0}

In this section, we generalize our findings to spline functions of
order $m \ge 1$.  Let us recall that the B-splines $B_{j,m}$ in
  \eqref{splines} being generated by the knot sequence
  $T_1 < \cdots < T_{N+m}$ are recursively defined by
\begin{equation*}
  B_{j,m}\mleft(t\mright) 
  \coloneqq \tfrac{t-T_j}{T_{j+m-1}-T_j} \, B_{j,m-1} \mleft(t\mright) 
  + \tfrac{T_{j+m}-t}{T_{j+m}-T_{j+1}} \, B_{j+1,m-1} \mleft(t\mright)
\end{equation*}
with
\begin{equation*}
  B_{j,1} \mleft(t\mright)
  \coloneqq \Ind_{[ T_j, T_{j+1})} \mleft(t\mright)
  \coloneqq 
  \begin{cases}
    1 \quad & t \in [ T_j, T_{j+1}),\\
    0 \quad & \text{else,}
  \end{cases}
  \addmathskip
\end{equation*}
see for instance \cite[\p~131]{Boo78}.  Further, we notice that for
$0 \le k \le m-2$ the $k$th derivative of the spline $f$ in
\eqref{splines} is given by
\begin{equation}
  \label{eq:k-deriv-spline}
  \frac{\diff^k}{\diff{t^k}} \, f \mleft(t\mright)
  = \sum_{j=1}^{N+k} c_j^{(m-k)} \, B_{j,m-k} \mleft(t\mright)
  ,
  \addmathskip
\end{equation}
where the coefficients $c_j^{(m-k)}$ are recursively defined by
\begin{equation*}
  c_j^{(m-k)} 
  \coloneqq \left(m-k\right) \, 
  \frac{c_j^{(m-k+1)} - c_{j-1}^{(m-k+1)}}{T_{j+m-k}-T_j}
  \qquad (j=1, \ldots , N+k),
\end{equation*}
with the convention that   $c_{0}^{(m-k+1)} = c_{N+k}^{(m-k+1)} =0$, see \cite[\p~139]{Boo78}.  For $k=m-1$, equation
\eqref{eq:k-deriv-spline} coincides with a step function, i.e., with the right derivative of the
linear spline $f^{(m-2)}$.  Further, in a distributional manner, the
$m$th derivative of $f$ is given by
\begin{equation}
  \label{eq:m-deriv-spline}
  \frac{\diff^m}{\diff{t^m}} \, f \mleft(t\mright)
  = \sum_{j=1}^{N+m} c_j^{(0)} \, \delta (t-T_j)
\end{equation}
with the coefficients
\begin{equation*}
 c_{1}^{(0)} := c_{1}^{(1)}, \quad c_{N+m}^{(0)} := -c_{N+m-1}^{(1)}, \quad  c_j^{(0)} 
  \coloneqq c_j^{(1)} - c_{j-1}^{(1)} \qquad (j=2, \ldots , N+m-1),
  \addmathskip
\end{equation*}
and the Dirac delta distribution $\delta$.

Applying the Fourier transform to \eqref{eq:m-deriv-spline}, we now obtain
\begin{equation}
  \label{eq:Fou-m-deriv}
  \widehat{f^{\left(m\right)}} \mleft( \omega \mright)
  = \left(\I \omega \right)^m \widehat f \mleft( \omega \mright)
  = \sum_{j=1}^{N+m} c_j^{(0)} \e^{- \I \omega T_j}.
  \submathskip
\end{equation}
and thus
\begin{equation}
  \label{eq:Fou-int-m-deriv}
  \omega^{2m} \abs{\widehat f \mleft( \omega \mright)}^2
  = \sum_{j=1}^{N+m} \sum_{k=1}^{N+m} c_j^{(0)} \overline c_k^{(0)} \, 
  \e^{-\I \omega (T_j - T_k)}.
  \addmathskip
\end{equation}
Since the exponential sum on the right-hand side of
\eqref{eq:Fou-int-m-deriv} has exactly the same structure as the
exponential sum in \eqref{ex1}, we can immediately generalize
Theorem~\ref{the:uni-step-fct} by considering
\begin{equation}\label{pneu}
P(\omega) := \omega^{2m}\abs{\widehat f \mleft(\omega\mright)}^2 
  = \sum_{\ell=-\nicefrac{(N+m)(N+m-1)}{2}}^{\nicefrac{(N+m)(N+m-1)}{2}}  \gamma_{\ell} \, \e^{-\I \omega \tau_{\ell}}.
\end{equation}

\begin{Theorem}
  \label{the:uni-spline-fct}
  Let $f$ be a spline function of the form \eqref{splines} of order
  $m$, whose knot distances $T_j - T_k$ differ pairwise for
  $j,k \in \{1,\dots,N+m\}$ with $j \ne k$, and whose coefficients
  satisfy
  \raisebox{0pt}[0pt][0pt]{$\absn{c^{(0)}_1} \ne
    \absn{c^{(0)}_{N+m}}$}.
  Further, let $h$ be a step size such that
  $h (T_j - T_k) \in (-\pi, \pi)$ for all $j, k$.  Then $f$ can be
  uniquely recovered from its Fourier intensities
  $\absn{\Fourier [f](h\ell)}$ with
  $\ell = 0, \dots, \nicefrac{3}{2}(N+m)(N+m-1)$ up to trivial
  ambiguities.
\end{Theorem}

\begin{Proof}
  The statement can be established by proceeding in the same manner as
  in Section~\ref{sec:retr-step-funct}.  First we apply Prony's method
  to the given samples $(h\ell)^{2m} \absn{\Fourier [f] (h\ell)}^{2}$
  with $\ell = 0, \dots, \nicefrac{3}{2}(N+m)(N+m-1)$ in order to
  determine the coefficients and frequencies of $P(\omega)$ in
  \eqref{pneu}. In a second step, the values $c_j^{(0)}$ and $T_j$ in
  (\ref{eq:Fou-m-deriv}) can be determined analytically as discussed
  in Theorem~\ref{the:uni-step-fct}.  Reversing the definition of
  \raisebox{0pt}[0pt][0pt]{$c_j^{(m-k)}$}, we can finally compute the
  unknown coefficients \raisebox{0pt}[0pt][0pt]{$c_j^{(m)}$} by
    \begin{align*}
      c_j^{(1)} &= c_j^{(0)} + c_{j-1}^{(1)}
      &&(j=1, \ldots , N+m-1) \\
      \shortintertext{and}
      c_j^{(m-k+1)} &= \tfrac{T_{j+m-k}-T_j}{m-k} \, c_j^{(m-k)} +
                      c_{j-1}^{(m-k+1)}
      &&(j=1, \ldots , N+k-1)
    \end{align*}
     with $c_0^{(1)} \coloneqq 0$ and $c_0^{(m-k+1)} \coloneqq 0$,
  which finishes the proof.
  \qed
\end{Proof}

\begin{Corollary}
  \label{cor:uni-spline-fct}
  Almost all spline functions $f$ of order $m$ in
  \eqref{splines} can be uniquely recovered from their Fourier
  intensities $\abs{\Fourier [f]}$ up to trivial ambiguities.
\end{Corollary}

\section{Numerical experiments}
\setcounter{equation}{0}

Since the proofs of Theorem~\ref{the:uni-step-fct} and
Theorem~\ref{the:uni-spline-fct} are constructive, they can be
straightforwardly transferred to numerical algorithms to recover a
spline function from its Fourier intensity.  However, the classical
Prony method introduced in Section~\ref{sec:pronys-method} is
numerically unstable with respect to inexact measurements and to
frequencies lying close together.  For this reason, there are numerous
approaches to improve the classical method.  In order to verify
Theorem~\ref{the:uni-step-fct} and Theorem~\ref{the:uni-spline-fct}
numerically, we apply the so-called approximate Prony method (APM)
proposed by Potts and Tasche in \cite[Algorithm~4.7]{PT10} for
recovery of parameters of an exponential sum of the form
\begin{equation}
  \label{poly}
  P(\omega) = \sum_{\ell=-M}^{M} \gamma_{\ell} \, \e^{-\I \omega \tau_{\ell} } 
\end{equation}
with $\tau_{\ell} = - \tau_{-\ell}$ and $\gamma_{\ell} = \overline{\gamma}_{-\ell}$.
The algorithm can be summarized as follows, where
  the exact number $2M+1$
  of the occurring frequencies in (\ref{poly}) needs not be known
  beforehand.

\begin{Algorithm}[Approximate Prony method \cite{PT10}]
  \label{alg:app-prony-meth}
  \emph{Input:} upper bound $L\in \N$ of the number $2M+1$ of
  exponentials; measurements $P(h k)$ with $k=0, \dots, 2\breve
  M$ and $\breve M \ge L$; accuracies $\varepsilon_1$,
  $\varepsilon_2$, and $\varepsilon_3$.
  \begin{enumerate}
  \item Compute a right singular vector
    $\Vek \lambda^{(1)} \coloneqq (\lambda_k^{(1)})_{k=0}^L$
    corresponding to the smallest singular value of the rectangular
    Hankel matrix
    $\Mat H \coloneqq (P(h(k+m)))_{k,m=0}^{2N-L,L}$.
  \item Evaluate the roots
    \raisebox{0pt}[0pt][0pt]{$ z_{j}^{(1)} =r_j^{(1)} \, \e^{\I
        \omega_j^{(1)}}$}
    of the polynomial
    $\Lambda^{(1)} (z) \coloneqq \sum_{k=0}^L \lambda^{(1)}_k \, z^k$
    with $\omega^{(1)}_j \in [0, \pi)$ and
    \raisebox{0pt}[0pt][0pt]{$\absn{r_j^{(1)}-1} \le
      \varepsilon_1$}.
  \item Compute a right singular vector
    $\Vek \lambda^{(2)} \coloneqq (\lambda_k^{(2)})_{k=0}^L$
    corresponding to the second smallest singular value of the
    rectangular Hankel matrix
    $\Mat H \coloneqq (P(h(k+m)))_{k,m=0}^{2N-L,L}$.
  \item Evaluate the roots
    \raisebox{0pt}[0pt][0pt]{$z_{j}^{(2)} = r_j^{(2)} \, \e^{\I
        \omega_j^{(2)}}$}
    of the polynomial
    $\Lambda^{(2)} (z) \coloneqq \sum_{k=0}^L \lambda^{(2)}_k \, z^k$
    with $\omega^{(2)}_j \in [0, \pi)$ and
    \raisebox{0pt}[0pt][0pt]{$\absn{r_j^{(2)}-1} \le
      \varepsilon_1$}.
  \item Determine all frequencies of the form
    $\omega_\ell \coloneqq \nicefrac{1}{2} \, (\omega^{(1)}_j +
    \omega^{(2)}_k)$
    if there exist indices $j$ and $k$ with
    \raisebox{0pt}[0pt][0pt]{$\absn{\omega^{(1)}_j - \omega^{(2)}_k}
      \le \varepsilon_2$}, and denote the number of found frequencies  by $\widetilde{M}$.
  \item Compute the coefficients $\gamma_{\ell}$ as least squares solution
    of the over-determined linear system
    \begin{equation*}
      \sum_{\ell=-\widetilde M}^{\widetilde M} \gamma_{\ell} \, \e^{\I h k
        \tau_{\ell}}
      = P(h k)
      \qquad(k = 0 , \dots, 2\breve M)
      \addmathskip
    \end{equation*}
    with $\tau_{\ell} = -\tau_{-\ell} = \nicefrac{\omega_{\ell}}{h}$ by using the
    diagonal preconditioner
    \begin{equation*}
      \Mat D 
      \coloneqq \diag \left(\tfrac{1-\absn{k}}{\widetilde M
        +1}\right)_{k=-\widetilde M}^{\widetilde M}.
    \end{equation*}
  \item Delete all pairs $(\tau_{\ell}, \gamma_{\ell})$ with $\absn{\gamma_{\ell}}
    \le \varepsilon_3$.
  \item Repeat step 6 with respect to the remaining frequencies
    $\tau_{\ell}$.
  \end{enumerate}
  \emph{Output:} coefficients $\gamma_{\ell}$ and frequencies $\tau_{\ell}$.
\end{Algorithm}

A second adaption of the proof of Theorem~\ref{the:uni-spline-fct}
concerns the reconstruction of the coefficients $c_j^{(m)}$ from the
recovered coefficients $c_j^{(0)}$.  In order to describe the relation
between the coefficients as linear equation system, we define the
rectangular matrices $\Mat C^{(m-k)} \in \R^{(N+k-1) \times (N+k)}$
for $k=0, \ldots , m-1$ elementwise by
\begin{equation*}
  C_{j\ell}^{(m-k)} \coloneqq 
  \begin{cases}
    \tfrac{m-k}{T_{j+m-k}-T_j} & \ell=j ,\\[\fskip]
    \tfrac{k-m}{T_{j+m-k}-T_j} & \ell=j-1 ,\\
    \hfill 0 \hfill & \text{else,}
  \end{cases}
  \qquad\text{and}\qquad
  C_{j\ell}^{(0)} \coloneqq 
  \begin{cases}
    \phantom{-} 1 & \ell=j, \\
    -1 & \ell=j-1, \\
    \phantom{-} 0 & \text{else.}
  \end{cases}
\end{equation*}
Then, the recursion between the coefficients $c_j^{(m-k+1)}$ and
$c_j^{(m-k)}$ can be stated as
\begin{equation*}
  \Mat C^{(m-k)} \, \Vek c^{(m-k+1)} = \Vek c^{(m-k)},
\end{equation*}
where we use the coefficient vectors
$\Vek c^{(m-k)} \coloneqq (c_j^{(m-k)})_{j=1}^{N+k}$.  Instead of
computing the coefficients stepwise from left to right, we can
determine the coefficients \raisebox{0pt}[0pt][0pt]{$c_j^{(m)}$} by
solving the over-determined linear equation system
\begin{equation}
  \label{eq:equa-sys-coeff}
  \Mat C^{(0)} \cdots \Mat C^{(m-1)} \, \Vek c^{(m)} = \Vek c^{(0)}.
\end{equation}
With these modifications, we recover a spline function of order $m$
from its Fourier intensity by the following algorithm.

\begin{Algorithm}[Phase retrieval]
  \label{alg:phase-retri}
  \emph{Input:} Fourier intensities $\absn{\Fourier[f](h k)}$ with
  $k= 0,\dots,2 \breve M$, step size $h>0$, order $m \ge 0$ of the spline
  function, upper bound $L$ of the number $N+m$ of knots with
  $L(L-1)<\breve M$, accuracy $\varepsilon$.
  \begin{enumerate}
  \item Compute the squared Fourier intensity of the
    $m$th derivative of the spline at the given points by
    \begin{equation*}
       \absn{\Fourier [f^{(m)}](h k)}^{2} = (h k)^{2m} \absn{\Fourier
        [f](h k)}^{2} 
      \qquad
      (k = 0, \dots, 2 \breve M).
    \end{equation*}
  \item Apply the approximate Prony method
    (Theorem~\ref{alg:app-prony-meth}) to determine the knot distances
    $\tau_\ell$ with
    $\ell = -\nicefrac{(N+m)(N+m-1)}{2}, \dots,
    \nicefrac{(N+m)(N+m-1)}{2}$
    in increasing order and the corresponding coefficients
    $\gamma_\ell$.
  \item Update the reconstructed distances and coefficients by
    \begin{equation*}
      \tau_\ell \coloneqq \frac{\tau_\ell - \tau_{-\ell}}{2}
      \qquad\text{and}\qquad
      \gamma_\ell \coloneqq \frac{\gamma_\ell +
        \overline{\gamma}_{-\ell}}{2}
    \end{equation*}
    for $\ell = 0, \dots, \nicefrac{(N+m)(N+m-1)}{2}$. 
  \item Set $T_1 \coloneqq 0$,
    $T_{N+m} \coloneqq \tau_{\nicefrac{(N+m)(N+m-1)}{2}}$, and
    $T_{N+m-1} \coloneqq \tau_{(\nicefrac{(N+m)(N+m-1)}{2})-1}$;
    find the index $\ell^*$ with
    $\absn{\tau_{\ell^*} - T_{N+m} + T_{M+m-1}} \le
    \varepsilon$; and compute the corresponding coefficients by
    \begin{equation*}
      c_1^{(0)} \coloneqq \abs{\frac{\gamma_{\nicefrac{(N+m)(N+m-1)}{2}} \, 
      \overline{\gamma}_{(\nicefrac{(N+m)(N+m-1)}{2})-1)}}{\gamma_{\ell^*}}}^\frac{1}{2}
    \end{equation*}
    as well as
    \begin{equation*}
      c_{N+m}^{(0)} \coloneqq \frac{\gamma_{\nicefrac{(N+m)(N+m-1)}{2}}}{\overline c_1^{(0)}}
      \qquad\text{and}\qquad
      c_N^{(0)} \coloneqq \frac{ \gamma_{(\nicefrac{(N+m)(N+m-1)}{2})-1}}{\overline c_1^{(0)}}.
    \end{equation*}
    Initialize the lists of recovered knots and coefficients by
    \begin{equation*}
      T \coloneqq [T_1, T_{N+m}, T_{N+m-1}]
      \qquad\text{and}\qquad
      C^{(0)} \coloneqq [c_1^{(0)}, c_{N+m}^{(0)}, c_{N+m-1}^{(0)}], 
    \end{equation*}
    and remove the used knot distances from the set ${\cal T} :=\{\tau_\ell\}_{\ell=0}^{\nicefrac{(N+m)(N+m-1)}{2}}$.
  \item For the maximal remaining distance $\tau_{k^*}$ in
    ${\cal T}$, determine the index $\ell^*$ with
    $\absn{\tau_{k^*} + \tau_{\ell^*} - T_{M+n}} \le
    \varepsilon$.
    \begin{enumerate}
    \item If $\absn{\tau_{k^*} - \tau_{\ell^*}} \le
      \varepsilon$, the knot distance corresponds to the centre
      of the interval $[T_1, T_{M+n}]$.  Thus append $T$ by
      $\nicefrac{T_{N+m}}{2}$ and $C^{(0)}$ by $\nicefrac{\gamma_{k^*}}{\overline c_1^{(0)}}$.
    \item Otherwise, compute the values
      $d^{(\mathrm r)} \coloneqq \nicefrac{\gamma_{k^*}}{\overline c_1^{(0)}}$
      and
      $d^{(\mathrm l)} \coloneqq \nicefrac{\gamma_{\ell^*}}{\overline c_1^{(0)}}$.  If
      \begin{equation*}
        \abs{c_{N+m}^{(0)} \, \overline d^{(\mathrm r)} - \gamma_{\ell^*}} < \abs{c_{N+m}^{(0)} \, \overline
          d^{(\mathrm l)} - \gamma_{k^*}},
    \end{equation*}
    then assume that \eqref{eq:case2} with $d^{(\mathrm r)}$,
      $\gamma_{k^*}$, $c_{N+m}^{(0)}$ instead of $c_{N-2}^{(0)}$,
      $\gamma_{(\nicefrac{N(N-1)}{2})-2}$, $c_N^{(0)}$ holds true and
    append $T$ by
    $\nicefrac{1}{2} \, (\tau_{k^*} + T_{N+m} - \tau_{\ell^*})$ and
    $C^{(0)}$ by $d^{(\mathrm r)}$, else assume that
    \eqref{eq:case1} with $d^{(\mathrm l)}$, $\gamma_{k^*}$,
      $c_{N+m}^{(0)}$ instead of $c_{2}^{(0)}$,
      $\gamma_{(\nicefrac{N(N-1)}{2})-2}$, $c_N^{(0)}$ holds true and
    append $T$ by
    $\nicefrac{1}{2} \, (\tau_{\ell^*} + T_{N+m} - \tau_{k^*})$ and
    $C^{(0)}$ by $d^{(\mathrm l)}$.
    \end{enumerate}
    Remove all distances between the new knot and the already
    recovered knots from ${\cal T}$ and repeat step 5 until
    the set ${\cal T}$ is empty.
  \item Determine the coefficients $c^{(m)}_j$ by solving the
    over-determined equation system \eqref{eq:equa-sys-coeff}.
  \end{enumerate}
  \emph{Output:} knots $T_j$ and coefficients $c_j^{(m)}$ of the
  signal \eqref{spikes} $(m=0)$ or the spline function in \eqref{splines} $(m>0)$.
\end{Algorithm}

\begin{table}[p]
  \centering
  \small
    \begin{tabular}{|rrr|rrr|rrr|}
      \hline
      \multicolumn{1}{|c}{$j$} 
      &\multicolumn{1}{c}{$T_j$} 
      &\multicolumn{1}{c}{$c_j^{(0)}$} 
      &\multicolumn{1}{|c}{$j$} 
      &\multicolumn{1}{c}{$T_j$} 
      &\multicolumn{1}{c}{$c_j^{(0)}$} 
      &\multicolumn{1}{|c}{$j$} 
      &\multicolumn{1}{c}{$T_j$} 
      &\multicolumn{1}{c|}{$c_j^{(0)}$ } 
      \\ \hline\hline
      1 & -53.5895 &  4.910 + 0.000i 
      & 6 &-28.1475& 0.278 + 0.598i
      & 11 & 1.3755& 2.887 + 3.828i
      \\
      2 & -50.2765 & -0.165 + 0.814i
      & 7 &-22.6005&-1.450 + 3.246i
      & 12 &20.0945&-1.423 + 0.397i 
      \\
      3 &-49.3765&-2.368 -- 1.314i
      & 8 &-19.6495&0.508 + 0.243i
      & 13 &33.4525& 0.023 -- 2.039i
      \\
      4 &-42.6915&-0.293 + 0.541i
      & 9 &-6.1705&0.073 -- 0.528i
      & 14 &34.8415&-2.997 + 3.767i 
      \\
      5 &-28.3915&-1.841 + 2.589i
      & 10 &-3.8985&3.135 + 0.339i
      & 15 & 53.5895&-0.064 -- 0.368i
      \\ \hline
    \end{tabular}
  \caption{Knots $T_j$ and coefficients $c_j^{(0)}$ of the spike
    function in Example~\ref{ex:spike-fun}}
  \label{tab:spike-fun}
\end{table}
\begin{figure}[p]
  \centering
  \begin{tabular}{@{}c@{\hspace{4pt}}c@{}}
    \subfloat[{Real part of the recovered and true spike function}]{
    \makebox[5.7cm]{\includegraphics%
    {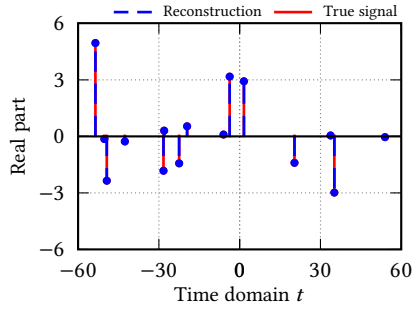}}
    }
    &
      \subfloat[{Imaginary part of the recovered and true spike function}]{
      \makebox[5.7cm]{\includegraphics%
      {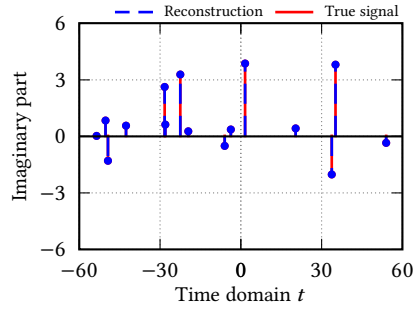}}
      }
    \\
    \subfloat[{Absolute error of the recovered knots}]{
    \makebox[5.7cm]{\includegraphics%
    {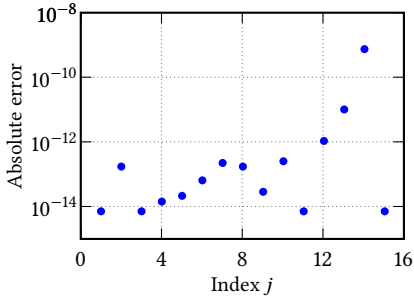}}
    }
    &
      \subfloat[{Absolute error of the recovered coefficients}]{
      \makebox[5.7cm]{\includegraphics%
      {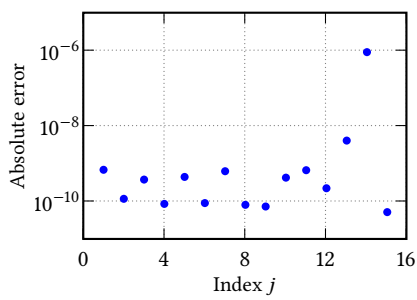}}
      }
  \end{tabular}
  \caption{Results of Algorithm~\ref{alg:phase-retri} for the spike
    function in Example~\ref{ex:spike-fun}}
  \label{fig:spike-fun}
\end{figure}

\begin{Example}
  \label{ex:spike-fun}
  In the first numerical example, we consider a spike function as
    in \eqref{spikes} with 15 spikes.  More precisely, the locations
    $T_j$ and the coefficients \raisebox{0pt}[0pt][0pt]{$c_j^{(0)}$}
    of the true spike function $f$ are given in
    Table~\ref{tab:spike-fun}.  In order to recover $f$ from the
    Fourier intensity measurements $\absn{\Fourier[f](h \ell)}$ with
    $\ell=0, \dots, 1000$ and with
    $h \approx 3.655 \, 073 \cdot 10^{-2}$, we apply
    Algorithm~\ref{alg:phase-retri} with the accuracies
    $\varepsilon \coloneqq 10^{-3}$,
    $\varepsilon_1 \coloneqq 10^{-5}$,
    $\varepsilon_2 \coloneqq 10^{-7}$, and
    $\varepsilon_3 \coloneqq 10^{-10}$.  The results of the phase
    retrieval algorithm and the absolute errors of the knots and
    coefficients of the recovered spike function are shown in
    Figure~\ref{fig:spike-fun}.  Although the approximate Prony method
    has to recover 211 knot differences, the knots and coefficients of
    $f$ are reconstructed very accurately.
  \qed
\end{Example}

\begin{figure}[p]
  \centering
  \begin{tabular}{@{}c@{\hspace{4pt}}c@{}}
    \subfloat[{Real part of the recovered and true spline function}]{
    \makebox[5.7cm]{\includegraphics%
    {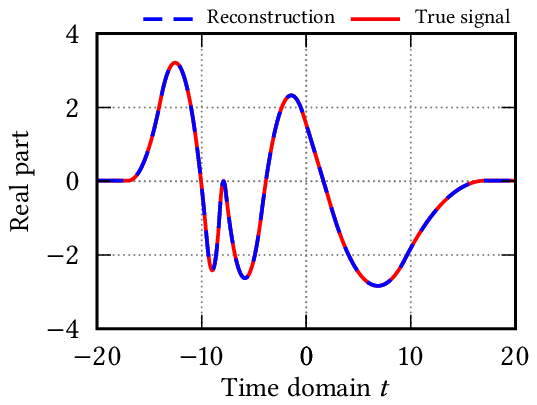}}
    }
    &
      \subfloat[{Imaginary part of the recovered and true spline
      function}]{
      \makebox[5.7cm]{\includegraphics%
      {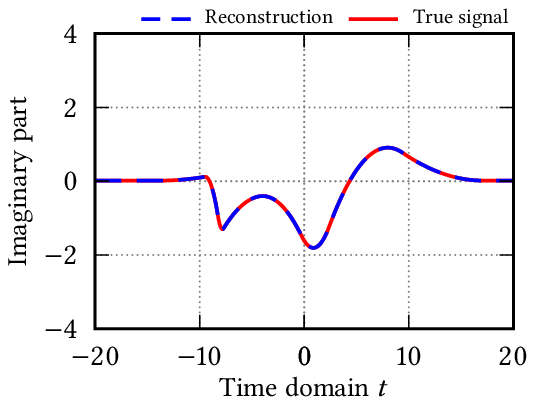}}
      }
    \\
    \subfloat[{Absolute error of the recovered knots}]{
    \makebox[5.7cm]{\includegraphics%
    {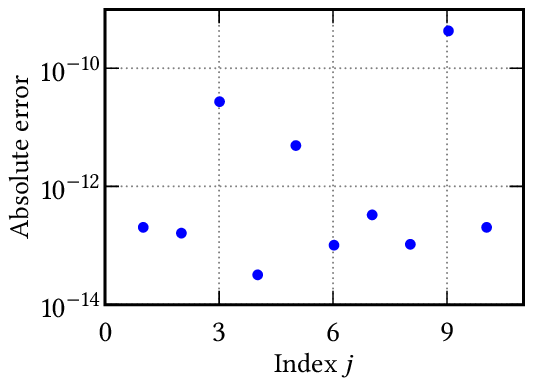}}
    }
    &
      \subfloat[{Absolute error of the recovered coefficients}]{
      \makebox[5.7cm]{\includegraphics%
      {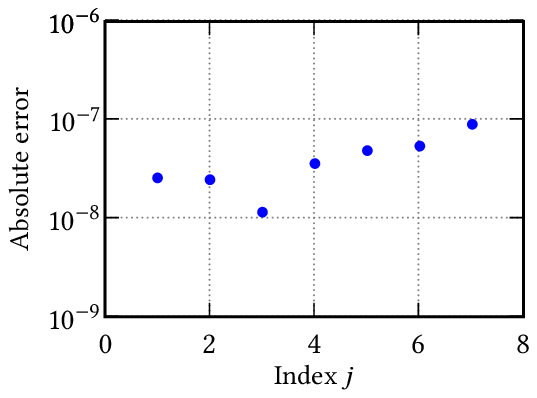}}
      }
  \end{tabular}
  \caption{Results of Algorithm~\ref{alg:phase-retri} for the spline
    function in Example~\ref{ex:spline-fun}}
  \label{fig:spline-fun}
\end{figure}

\begin{Example}
  \label{ex:spline-fun}
 In the second example, we apply Algorithm~\ref{alg:phase-retri}
    to recover the piecewise quadratic spline function ($m=3$) in
    \eqref{splines} with the knots and coefficients in
    Table~\ref{tab:spline-fun} from the Fourier intensity measurements
    $\absn{\Fourier[f](h \ell)}$ with $\ell=0, \dots, 400$ and with
    $h \approx 3.088 \, 663 \cdot 10^{-2}$.  As accuracies for the
    phase retrieval algorithm and the approximate Prony method, we
    choose $\varepsilon \coloneqq 10^{-3}$,
    $\varepsilon_1 \coloneqq 10^{-5}$,
    $\varepsilon_2 \coloneqq 10^{-10}$, and
    $\varepsilon_3 \coloneqq 10^{-10}$.  In
    Figure~\ref{fig:spline-fun}, the recovered function is compared
    with the true signal.  Again, the reconstructed knots and
    coefficients have only very small absolute errors.
  \qed
\end{Example}

\begin{table}[p]
  \centering
  \small
    \begin{tabular}{|rrr|rrr|rrc|}
      \hline
      \multicolumn{1}{|c}{$j$} 
      &\multicolumn{1}{c}{$T_j$} 
      &\multicolumn{1}{c}{$c_j^{(3)}$} 
      &\multicolumn{1}{|c}{$j$} 
      &\multicolumn{1}{c}{$T_j$} 
      &\multicolumn{1}{c}{$c_j^{(3)}$} 
      &\multicolumn{1}{|c}{$j$} 
      &\multicolumn{1}{c}{$T_j$} 
      &\multicolumn{1}{c|}{$c_j^{(3)}$ } 
      \\ \hline\hline
      1 & -17.022 & 5.342 + 0.000i
      & 5 &-7.745&3.597 -- 0.334i
      & 8 &2.309& ---
      \\
      2 &-13.921 &  -3.569 + 0.132i
      & 6 &-4.313&0.554 -- 2.251i
      & 9 &9.318& ---
      \\
      3 &-9.536 &  0.440 -- 1.413i
      & 7 &-0.336&-4.072 + 1.433i
      & 10 & 17.022& ---
      \\
      4 & -8.301&  -4.685 -- 0.499i
      &&&
      &&& \phantom{3.5970 -- 0.334i}
      \\ \hline
    \end{tabular}
  \caption{Knots $T_j$ and coefficients $c_j^{(3)}$ of the spline
    function in Example~\ref{ex:spline-fun}}
  \label{tab:spline-fun}
\end{table}

\section*{Acknowledgements}
The first author gratefully acknowledges the funding of this work by
the Austrian Science Fund (FWF) within the project P~28858, and the
second author the funding by the German Research Foundation (DFG)
within the project PL~170/16-1.  The Institute of Mathematics and
Scientific Computing of the University of Graz, with which the first
author is affiliated, is a member of NAWI Graz
(\texttt{http://www.nawigraz.at/}).

\bibliographystyle{alphadinUK}
{\footnotesize \bibliography{LITERATURE}}

\end{document}